\documentclass[preprint]{elsarticle}
\usepackage[matrix,arrow,curve]{xy}
\usepackage{graphicx}
\usepackage{amsmath}
\usepackage{amssymb}
\usepackage[utf8]{inputenc}
 \usepackage[T2A]{fontenc}
\usepackage[english]{babel}
\usepackage{bm}
\usepackage{color}
\usepackage{ulem}
\usepackage{xcolor}
\newtheorem{thm}{Theorem}
\newtheorem{lem}{Lemma}

\newproof{pf}{Proof}

\newdefinition{rmk}{Remark}

\begin{document}

\begin{frontmatter}
\title{Differential invariants of feedback transformations for quasi-harmonic oscillation equations}

\author[1,2]{D.S.~Gritsenko}\ead{gricenko@physics.msu.ru}
\address[1]{Faculty of Physics, Lomonosov Moscow State University, GSP-2, Leninskie Gory, Moscow 119991, Russia}
\address[2]{V. A. Trapeznikov Institute of Control Sciences of Russian Academy of Sciences, 65 Profsoyuznaya street, Moscow 117997, Russia}
\author[3,1]{O.M.~Kiriukhin}\ead{kiriukhin@uchicago.edu}
\address[3]{University of Chicago Booth School of Business, 5807 South Woodlawn Avenue Chicago, Illinois 60637, USA}
\journal{Journal of Differential Equations}

\begin{abstract}
The goal and the main result of the paper is to provide a complete description of the field of rational differential invariants of one class of second order ordinary differential equations with scalar control parameter with respect to Lie pseudo-group of local feedback transformations. In particular, considered class describes behavior of conservative mechanical systems. We construct the class of rational differential invariants that separate regular orbits. It is well known that differential invariants form algebra with respect to the operation of addition and multiplication \cite{Al1988}. In our case, constructed rational differential operators form a field (in algebraic sense). Rational differential invariants were studied by M. Rosenlicht \cite{Rose1956, Rose1963}, B. Kruglikov and V. Lychagin \cite{Kruglikov}.

\end{abstract}

\end{frontmatter}

\section{INTRODUCTION}
\par The feedback classification and equivalence problem has drawn  attention of the researchers for several decades. A variety of approaches for its treatment has been proposed, developed and documented in extensive existing literature. They are the following. The first one is Cartan's equivalence method originally proposed by Ellie Cartan in \cite{Cartan} and developed by his numerous heirs (\cite{G,Ku}). Another one is based on Hamiltonian formalism for optimal control systems (See \cite{B,J2}). One more is the geometric approach that uses the equivalence of the distributions and affine distributions. It was first formulated by Jacubczyk in \cite{J1} and followed by \cite{Br}--\cite{ZR}. An alternative approach is the formal one proposed by Kang in \cite{KK} and further studied in \cite{K}--\cite{TR}. In this paper we employ the approach based on differential invariants introduced in \cite{VKL1986} and developed in \cite{PO}.

\par Our paper contributes to the existing literature in several ways. We provide a complete description of the field of Petrov invariants  of {\it quasi-harmonic oscillation equation} (QHOE) with respect to Lie pseudo-group of local feedback transformations. We also provide straightforward algorithm for calculating the Petrov invariants up to the $k$-th order. Furthermore, we specify some simple free coordinate forms of QHOE based on the derived differential invariants. 

\par Feedback transformations of the control-parameter-dependent systems are analogous to Lie transformations of the differential equations. They are often used for the reduction of such systems to canonical or normal  forms and simplification of equations.

\par The structure of this paper is as follows. In Section 2 we derive a set of basal vector fields  and local translation groups (local feedback transformations, local diffeomorphisms) along their trajectories. These transformations preserve a class of QHOE. Finally, the action of local feedback transformations on QHOE is given. In Section 3 we provide an algorithm for the calculation of the differential invariants up to the $k$-th order. Using this algorithm we find second- and third-order differential invariants of QHOE. Finally we specify some simple free coordinate forms of QHOE. In Section 4 we define G-invariant derivation and find their explicit form for QHOE. We also obtain syzygies for the third-order differential invariants. In Section 5 we show that differential invariants obtained in Section 3 are Petrov invariants. We then calculate the dimension of algebra of Petrov invariants. Finally, we formulate a theorem that provides a complete description of the field of Petrov invariants. In Section 6 we summarize main results of the paper and discuss further possible studies.

\par To perform the required calculations MAPLE 15 software has been employed, especially DifferentialGeometry and JetCalculus packages developed by I. Anderson.

 \section{ADDMISSIBLE FEEDBACK TRANSFORMATIONS}
\par It is generally known that the classification problem for the ordinary differential equations (ODEs) with respect to Lie transformations is among the key ones in the theory of the differential equations. Construction of algebra of differential invariants is the first step to classification of equations.
 \par  In this paper we construct the field of rational differential invariants with respect to the local feedback transformations
\begin{equation} \label{transformation}
 \varphi: (x,y,u)\longmapsto (X(x,y),Y(x,y),U(u))
 \end{equation}
for equations of the following form:
 \begin{equation}\label{mainEq}
 \frac{d^{2}y}{d x^{2}}+f(y,u)=0.
 \end{equation}

\par Here $y=y(x)$ and $u=u(x)$ are scalar functions of a scalar argument $x$. The function $f(y,u)$ is smooth, i.e. of class $C^\infty$. The function $u=u(x)$ is a control parameter. Equations (\ref{mainEq})  describe, for example, the mechanical one-dimensional conservative systems. We call equations  (\ref{mainEq})  {\it{control-parameter-dependent quasi-harmonic oscillaton equations}} (QHOE).

\par Having followed an approach proposed by Sophus Lie \cite{Lie2011}, consider one-parameter local group of transformations
\begin{equation}\label{IZtransformation}
 \varphi_t: (x,y,u)\longmapsto (X_t(x,y),Y_t(x,y),U_t(u)),
 \end{equation}
instead of generalized point transformations (\ref{transformation}). Here $t$ is a scalar parameter $t\in I$, where $I\subset \mathbb{R}$ is an open interval, $0\in I$ and $X_t(x,y),Y_t(x,y),U_t(u)$ are one-parameter families of smooth functions which depend on the parameter $t$ smoothly as well.

\par We suppose that at the point $t=0$ transformation (\ref{IZtransformation}) is identical, i.e. $ X_0(x,y)=x, \ \ Y_0(x,y)=y, \ \ U_0(u)=u$.
Transformations (\ref{IZtransformation}) are generated by vector fields of the form
\begin{equation}\label{VectorFieldX}
 X=A(x,y)\frac{\partial}{\partial x}+B(x,y)\frac{\partial}{\partial y}+C(u)\frac{\partial}{\partial u},
\end{equation}
i.e. the transformation $\varphi_t$ is a translation along trajectories of the vector field $X$ from $t=0$ to $t$. Transformations (\ref{IZtransformation}) form Lie pseudo-group.

Let us find  feedback transformations (\ref{transformation}) preserving a class of QHOE. Denote the space of 2-jets of smooth vector-functions $$\mathbb{R}\rightarrow \mathbb{R}^2, \quad x \mapsto (y(x),u(x))$$ by $J^2(1,2)$. Let $x, y, u, y_x, u_x, y_{xx}, u_{xx}$ be canonical coordinates on $J^2(1,2)$ (see \cite{VKL1986}). This means that the Cartan forms have the following coordinate representations:
$$\omega_1=dy-y_xdx,\, \omega_2=du-u_xdx,\, \omega_3=dy_x-y_{xx}dx,\, \omega_4=du_x-u_{xx}dx.$$
These 1-forms define the Cartan distribution $\mathcal{C}$ on $J^2(1,2)$:
$$
\mathcal{C}: J^2(1,2)\ni \theta \mapsto \mathcal{C}(\theta)=\bigcap_{i=1}^4\ker\omega_i\subset T_\theta J^2(1,2).
$$
Here $T_\theta J^2(1,2)$ is a tangent space of $J^2(1,2)$ at the point $\theta$.

\par Equation (\ref{mainEq}) defines a hyper-surface
\begin{equation*}
 \mathcal{E}_f=\{y_{xx}+f(y,u)=0\} \subset J^2(1,2).
\end{equation*}

\par Now let us  find the vector fields (\ref{VectorFieldX}), whose flows are of the form (\ref{IZtransformation}) and  preserve a class of equations  (\ref{mainEq}). Thereafter the prolongations of vector fields and infinitesimal transformations in the space of 2-jets are denoted by $X^{(2)}$ and  $\varphi_t^{(2)}$, respectively (See \cite{VKL1986}). Transformations (\ref{IZtransformation}) preserve the class of equations (\ref{mainEq}) if and only if

 \begin{equation}\label{equality}
 \left( \varphi_t^{(2)}\right)^*(y_{xx}+f(y,u))=\lambda_t(y_{xx}+g_t(y,u)),
 \end{equation}
where $\lambda_t$ is a local one parameter family of smooth functions on the space $J^2(1,2)$, and $g_t(y,u)$ is a local one parameter family of smooth functions of variables $y$ and  $u$, such that $\lambda_0=1$ and $g_0(y,u)=f(y,u)$. Taking the derivative of (\ref{equality}) with respect to $t$ at $t=0$, we obtain:

\begin{align}\label{ld}
& \left. \frac{d}{dt}\right|_{t=0} \left( \varphi^{(2)}_t \right)^* (y_{xx}+f(y, u)) = \lambda_0 \left. \frac{d}{dt}\right|_{t=0} (y_{xx}+g_t(y, u)) + \\ \nonumber
& +(y_{xx}+g_0(y, u))\left. \frac{d\lambda_t}{dt}\right|_{t=0}=\left. \frac{d}{dt}\right|_{t=0} g_t(y, u) + (y_{xx}+f(y, u))\left. \frac{d\lambda_t}{dt}\right|_{t=0}.\label{equal}
\end{align}

\par The left side of (\ref{ld}) is the Lie derivative along the vector field $X^{(2)}$ of  the function $y_{xx}+f(y_0,u)$. A restriction of the last equality of (\ref{ld}) on $\mathcal{E}_f$ is:
 \begin{equation*}
 \left. L_{X^{(2)}}(y_{xx}+f(y,u))\right|_{\mathcal{E}_f} = G(y,u),
 \end{equation*}
 or
  \begin{equation}\label{LieDerivative}
 \left.X^{(2)}(y_{xx}+f(y, u))\right|_{\mathcal{E}_f} = G(y, u),
 \end{equation}
 where

  \begin{equation*}
G(y,u)=\left.\frac{d}{dt}\right|_{t=0} g_t(y, u).
 \end{equation*}

Vector equation (\ref{LieDerivative}) can be rewritten as a system of scalar linear differential equations with respect to functions $A$, $B$:
\[
\left\{
\begin{array}{l}
A_{{yy}}=0,\\
B_{{yy}}-2\,A_{{xy}}=0,\\
2\,B_{{xy}}-A_{{xx}}-3\,fA_{{y}}=0,\\
B_{{xx}}-Cf_{{u}}-Bf_{{y}}-G \left( y,u \right) +fB_{{y}}-2\,fA_{{x}}=0.
 \end{array}
\right.
\]

These equations must be satisfied identically for  all functions $f$ and some function $G$. The general solution of this system is:

$$
A(x, y)=\alpha x+\beta, \quad B(x, y)=\gamma+\delta y.
$$
Here  $\alpha,\beta,\gamma,\delta$ are arbitrary constants. Hence, any  vector fields preserving the class of equations  (\ref{mainEq}) can be presented as a linear combination with constant coefficients of the following  basal vector fields:

 \begin{equation*}
 X_1=\frac{\partial}{\partial x},\quad X_2=\frac{\partial}{\partial y}, \quad X_3=y\frac{\partial}{\partial y}, \quad X_4=x\frac{\partial}{\partial x},\quad X_5=C(u)\frac{\partial}{\partial u}.
 \end{equation*}

 Local translation groups along these fields can be written as follows:
 \begin{align*}
 \varphi_{1, t}& :(x, y, u)\longmapsto( x+t, y, u), \quad &\varphi_{2, t}& :(x, y, u)\longmapsto( x, y+t, u),\\
 \varphi_{3, t}& :(x, y, u)\longmapsto( x, e^t y, u),\quad &\varphi_{4,t}& :(x, y, u)\longmapsto( e^t x, y, u),\\
 \varphi_{5, t}& :(x, y, u)\longmapsto( x, y, U_t(u)).
\end{align*}

Let us calculate how the transformations  $\varphi_{i,t}$ act on equation (\ref{mainEq}) and on the function $f$. Transformation $\varphi_{1,t}$  doesn't change the form of equation (\ref{mainEq}). Applying transformations $\varphi_{2,t}^{(2)},\dots, \varphi_{5,t}^{(2)}$ to the left side of (\ref{mainEq}) result in:
\begin{align*}
 \left( \varphi_{2, t}^{(2)}\right)^{*} \left( y_{xx}+f(y, u)\right) & =y_{xx}+f(y+t, u) \\
 \left( \varphi_{3, t}^{(2)}\right)^{*} \left( y_{xx}+f(y, u)\right) & =y_{xx}+e^{-t}f(y, u),\\
 \left( \varphi_{4, t}^{(2)}\right)^{*} \left( y_{xx}+f(y, u)\right) & =y_{xx}+e^{\frac{t}{2}}f(y, u),\\
 \left( \varphi_{5, t}^{(2)}\right)^{*} \left( y_{xx}+f(y, u)\right) & =y_{xx}+f(y, U_t(u)).
 \end{align*}

 \section{DIFFERENTIAL INVARIANTS}

 Construct a trivial vector bundle 

 \begin{equation*}
 \pi: \mathbb{R}^3\rightarrow \mathcal{B}=\mathbb{R}^{2}, \qquad \pi: (y,u,z)\mapsto (y,u).
 \end{equation*}

Sections
 \begin{equation*}
 s_f: (y, u) \mapsto (y, u, f(y,u))
 \end{equation*}
of this bundle are smooth functions on the base $\mathbb{R}^2$: $z=f(y,u)$. Transformations $\varphi_{2,t}^{(2)},\dots, \varphi_{5,t}^{(2)}$ form a Lie pseudo-group and act on the total space of $\pi$. This pseudo-group  we denote by $G$.  Corresponding vector fields are:
  \begin{equation*}
 Y_1=\frac{\partial}{\partial y},\quad   Y_2=y\frac{\partial}{\partial y}, \quad  Y_3=z\frac{\partial}{\partial z},\quad   Y_4=H(u)\frac{\partial}{\partial u}.
 \end{equation*}

\par Let  $J^k(\pi)$ be the space of $k$-jets of sections of the bundle  $\pi$ and $y, u, z, z_{y}, z_{u}, z_{yy}, z_{yu}, z_{uu},\ldots, z_{u\dots u}$
be the canonical coordinates  on  this space. Let $\varphi^{(k)}$ and $Y^{(k)}$ be the prolongations to the space $J^k(\pi)$ of a transformation $\varphi\in G$ and a vector field $Y\in \mathcal{G}$ respectively. A set
$$
\mathcal{O}^k(\theta)=\bigcup_{\varphi\in G}\varphi^{(k)}(\theta)
$$
is called an {\it orbit} of the point $\theta\in J^k(\pi)$. For a point $\theta\in J^k(\pi)$ define a tangent subspace
$$
\mathcal{H}^k(\theta)=\mathrm{span}\bigcup_{H\in C^\infty(\mathbb{R})}\left(H(u)\frac{\partial}{\partial u}\right)^{(k)}
$$
of $T_\theta J^k(\pi)$. A point  $\theta\in J^k(\pi)$ is called {\it regular} if the dimension of the tangent subspace
\begin{equation}\label{Zk}
Z^k(\theta)=\mathrm{span}\left(\mathbb{R}Y_{1,\theta}^{(k)},  \mathbb{R}Y_{2,\theta}^{(k)},\mathbb{R}Y_{3,\theta}^{(k)}, \mathcal{H}(\theta)\right)\subset T_\theta J^k(\pi)
\end{equation}
is maximal and {\it singular} otherwise. An orbit is called {\it{regular}} if all its points are regular.

Consider, for instance, the space of 0-jets $J^0(\pi)$. If a point $\theta$ belongs to the hyperplane $\{z=0\}\subset J^0(\pi)$ then $\dim Z^0(\theta)=2$ and $\dim Z^0(\theta)=3$ for other points. We see that any point of the hyperplane  $\{z=0\}$ is singular. Moreover, since the Lie group $G$ acts transitively on $\{z=0\}$,  this hyperplane is a singular orbit. Hyperplane $\{z=0\}$ divides the space $J^0(\pi)$ into two connected components. Lie pseudo-group $G$ operates transitively at any connected component. Points of the hyperplanes $\{z=0\}$, $\{z_y=0\}$ and $\{z_u=0\}$ are singular as well.

A function $J$  on the space $J^{k}(\pi)$ smooth in its domain of definition and rational with respect to variables $z_{x}, z_{y}, z_{xx}, z_{xy}, \ldots$ is called {\it{ differential invariant of order  $\leq k$}} of the Lie pseudo-group $G$ if it is a constant on the orbits of prolonged Lie pseudo-group $G^{(k)}$, or, equivalently,
$$
\left(\varphi^{(k)}\right)^\ast(J)=J.
$$ (See \cite{Al1988}).

Rational differential invariants form a field (in algebraic sense) \cite{KL2013}. Following \cite{KL2013} rational differential invariants of {\it QHOE} we call {\it{Petrov invariants}}.

Having solved  the system of differential equations
\begin{equation}
 Y^{(k)}(J) = 0, \label{invariants}
 \end{equation}
for any $Y \in\mathcal{G}$ we find differential invariants of order $\leq k$ of  $G$. First non-trivial differential invariants appear on the space of 2-jets. We find them solving system (\ref{invariants}) for $k=2$.  The prolongations of the vector fields $Y_i$ to the space of 2-jets can be written as:

\begin{align*}
 Y_1^{(2)} & = Y_1,\\
 Y_2^{(2)} & = Y_2 - z_{y} \frac{\partial}{\partial z_{y}} - 2 z_{yy} \frac{\partial}{\partial z_{yy}} - 2 z_{yu} \frac{\partial}{\partial z_{yu}}, \\
 Y_3^{(2)} & = Y_3 + z_{y} \frac{\partial}{\partial z_{y}} + z_{u} \frac{\partial}{\partial z_{u}} + z_{yy} \frac{\partial}{\partial z_{yy}} + z_{yu} \frac{\partial}{\partial z_{yu}} + z_{uu} \frac{\partial}{\partial z_{uu}}, \\
 Y_4^{(2)} & = Y_4 - H_{u}(u) z_{u} \frac{\partial}{\partial z_{u}} - H_{u}(u) z_{yu} \frac{\partial}{\partial z_{yu}} - \Bigl( H_{uu}(u) z_{u} +2 H_{u}(u) z_{uu} \Bigr) \frac{\partial}{\partial z_{uu}}.
 \end{align*}

The system (\ref{invariants}) being solved at $k=2$ gives two  second-order differential invariants:

  \begin{equation} \label{inv2}
 J_{21} =\frac{z_{yy}z}{z_y^2},\quad J_{22} =\frac{z_{yu}z}{z_{y}z_u}
 \end{equation}

\par The prolongations of vector fields $Y_i$ to the space of 3-jets  can be written as:

\begin{align*}
 Y_1^{(3)} & = Y_1, \\
 Y_2^{(3)} & = Y_2^{(2)} - 3 z_{yyy} \frac{\partial}{\partial z_{yyy}} - 2 z_{yyu}\frac{\partial}{\partial z_{yyu}} - z_{yuu}\frac{\partial}{\partial z_{yuu}}, \\
 Y_3^{(3)} & = Y_3^{(2)} + z_{yyy}\frac{\partial}{\partial z_{yyy}} + z_{yyu}\frac{\partial}{\partial z_{yyu}} + z_{yuu}\frac{\partial}{\partial z_{yuu}} + z_{uuu}\frac{\partial}{\partial z_{uuu}}, \\
 Y_4^{(3)} & = Y_4^{(2)} - H_{u}(u) z_{yyu} \frac{\partial}{\partial z_{yyu}} - \Bigl(H_{uu}(u) z_{yu} + 2 H_{u}(u) z_{yuu} \Bigr)\frac{\partial}{\partial z_{yuu}} - \\ & - \Bigl(H_{uuu}(u) z_{u} + 3 H_{uu}(u) z_{uu} + 3 H_{u}(u) z_{uuu} \Bigr) \frac{\partial}{\partial z_{uuu}} .
 \end{align*}

Resolving system (\ref{invariants})  for $k=3$ we obtain three third-order differential  invariants:

 \begin{equation*}\label{inv3}
 J_{31}  =\frac{z_{yyy}z^2}{z_y^3}, \quad J_{32}  =\frac{z_{yyu}z^2}{z_y^2z_u}, \quad J_{33}  =\frac{z_{yuu}z^2}{z_u^2z_y}-J_{22}\frac{z_{uu}z}{z_{u}^{2}}.
 \end{equation*}
Expressions for higher-order differential invariants are given in Appendix A. Constructed invariants allow us to describe some QHOE in free-coordinate forms.

 \begin{thm}
1. Equation (\ref{mainEq})  is locally equivalent to the equation
 \begin{equation*}
  \frac{d^2 y}{d x^2}+ (\alpha(u)y)^n = 0
  \end{equation*}
with respect to feedback transformations if  and only if
\begin{equation}\label{nf1}
\left\{
\begin{array}{ll}
 J_{21}(f)    &= \frac{n-1}{n},  \quad J_{22}(f)  = 1, \quad J_{31}(f)  = \frac{n^2 -3n +2}{n^2},\\
  J_{32}(f)   &= \frac{n-1}{n}, \quad J_{33}(f)  = 0
\end{array}
\right.
\end{equation}
 for some natural number $n$.

2.  Equation (\ref{mainEq})  is locally equivalent to the equation  \[ \frac{d^2 y}{d x^2}+\alpha(u) y+\beta(u)=0\] if and only if
\begin{equation}\label{nf2}
 J_{21}(f)=0.
 \end{equation}

3. Equation  (\ref{mainEq})  is locally equivalent to the  equation  \[ \frac{d^2 y}{d x^2}+\alpha(u) y^2+\beta(u) y +\gamma(u)=0\] if and only  if
\begin{equation}\label{nf3}
 J_{31}(f)=0.
 \end{equation}

Here $\alpha, \beta, \gamma$ are of class $C^\infty$ and $J(f)$ is restriction of the invariant $J$ to the function $f$.
\end{thm}

\begin{pf}
To prove the theorem it is sufficient to solve differential equations (\ref{nf1}), (\ref{nf2}), (\ref{nf3}).
\end{pf}

\section{INVARIANT DERIVATIONS}

Differential operator
 \begin{equation}\label{operator}
 \nabla = M\frac{d}{dy} + N\frac{d}{du}
 \end{equation}
is called $G$-{\it invariant derivation} if it commutes with every element of any prolongation of the Lie algebra  $\mathcal{G}$, where $M$ and $N$  are the functions on $J^\infty(\pi)$.

It means that the following diagram

 $$
\xymatrix{
{C^\infty(J^\infty{{(\pi)}})} \ar[rrr]^{\nabla} \ar[dd]_{X^{(\infty)}} &&& {C^\infty(J^\infty{{(\pi)}})} \ar[dd]^{X^{(\infty)}} \\
&&&\\
{C^\infty(J^\infty{{(\pi)}})} \ar[rrr]_{\nabla} &&& {C^\infty(J^\infty{{(\pi)}})}
}
$$
 is commutative for any vector field $Y^{*} \in \mathcal{G}^{(\infty)}$. Here

 \begin{align*}
\frac{d}{dy} = \frac{\partial}{\partial y} + z_{y} \frac{\partial}{\partial z} + z_{yy} \frac{\partial}{\partial z_{y}} + z_{yyy} \frac{\partial}{\partial z_{yy}}  + \ldots,\\
 \frac{d}{du} = \frac{\partial}{\partial u} + z_{u} \frac{\partial}{\partial z} + z_{uu} \frac{\partial}{\partial z_{u}} + z_{uuu} \frac{\partial}{\partial z_{uu}}  + \ldots
 \end{align*}
are the operators of total differentiation by the variables $y$ and $u$ respectively (See \cite{KLR2007}).

Invariant  derivation  let us construct new differential invariants based on known ones. Indeed, let $J$ be the differential invariant and $\nabla$ be the invariant  derivation. Then

 \begin{equation*}
 Y^\ast(\nabla(J))=\nabla(Y^\ast(J))=0
 \end{equation*}
for any vector field $Y^\ast \in \mathcal{G}^{(\infty)}$. Therefore the function $\nabla (J)$ is  also differential invariant.

 \begin{lem} Differential operators
 \begin{equation}\label{invdiff}
 \nabla_1 = \frac{z}{z_{y}} \frac{d}{dy},\qquad\mbox{and}\qquad \nabla_2 = \frac{z}{z_{u}} \frac{d}{du}
 \end{equation}
 are $G$-invariant  derivations.
 \end{lem}

 \begin{pf}
 According to \cite{Lych2009} if the functions $M$ and $N$ satisfy the following system of differential equations:

 \[
\left\{
\begin{array}{ll}
 X^{(\infty)}(M) + \frac{d}{dy}\left(\frac{\partial h}{\partial z_y} \right)M + \frac{d}{du}\left(\frac{\partial h}{\partial z_u} \right)N& = 0,\\
 X^{(\infty)}(N) + \frac{d}{dy}\left(\frac{\partial h}{\partial z_y} \right)M + \frac{d}{du}\left(\frac{\partial h}{\partial z_u} \right)N & = 0,
\end{array}
\right.
\]
for any vector field $X \in \mathcal{G}$, then operator (\ref{operator}) is invariant  derivation. Here $h$ is a generating function of the vector field $X^{(1)}$, i.e. $h=\omega(X^{(1)})$, where
$$\omega= dz-z_ydy-z_udu$$ is the Cartan form on $J^1(\pi)$. Having resolved this system  restricted on the space of 2-jets for the vector fields  $Y_1,\ldots,Y_4$ we obtain:
$$
M = \frac{C_1 z}{z_{y}},\qquad N = \frac{C_2 z}{z_{u}},
$$
where $C_1, C_2$ are arbitrary constants. Having assumed,
 $C_1 = 1$, $C_2 = 1$ we obtain invariant   derivations (\ref{invdiff}).
 \end{pf}
\par Invariant  derivation is determined up to multiplication by a differential invariant. Note that
\begin{equation*}
 [\nabla_1,\nabla_2] = (-1 + J_{22}) \nabla_1 + (1 - J_{22}) \nabla_2,
 \end{equation*}
where $J_{22}$ is the second order differential invariant (see (\ref{inv2})). Applying the constructed invariant  derivations  $\nabla_{1,2}$ to the invariants (\ref{inv2}) we obtain:

 \begin{equation}\label{NablaOnInvariants}
 \left\{
  \begin{array}{ll}
 \nabla_1( J_{21}) & = J_{21} - 2 J_{21}^2 + J_{31}, \\
 \nabla_2( J_{21}) & =  J_{21} - 2 J_{21} J_{22} + J_{32}, \\
 \nabla_1 (J_{22}) & =  J_{22} - J_{21} J_{22}  - J_{22}^2 + J_{32}, \\
 \nabla_2 (J_{22}) & =  J_{22} - J_{22}^2 + J_{33}.
 \end{array}
 \right.
 \end{equation}
 
 Syzygies for higher-order differential invariants are given in Appendix B.

\section{STRUCTURE OF THE FIELDS OF INVARIANTS}

We say that a set of differential invariants $J_1,\dots, J_\nu$ of order $\leq k$ form a {\it{local basis}} of the field of Petrov invariants in an open domain $\mathcal{D}\subset J^k(\pi)$ if the following conditions hold:
\begin{enumerate}
\item the invariants $J_1,\dots, J_\nu$ are functionally independent in the domain $\mathcal{D}$, i.e. $dJ_1\wedge\dots\wedge dJ_\nu\neq 0$ in $\mathcal{D}$;
\item in the domain $\mathcal{D}$ any differential invariant of order $\leq k$ is a function of $J_1,\dots, J_\nu$ .
\end{enumerate}

In this case we say that the Petrov invariants $J_1,\dots, J_\nu$ are {\it{basic}} invariants of order $\leq k$ in the domain $\mathcal{D}$.

\begin{lem}\label{LemmaNu}
The number of basic Petrov invariants of order $\leq k$ $(k\geq 1)$ is
\begin{equation}\label{nu}
\nu(k)=\frac{k^2+k-2}{2}.
\end{equation}
\end{lem}
\begin{pf}
Construct projections
$$
\pi_{k,0}: J^k(\pi)\rightarrow \mathcal{B},\quad \pi_{k,r}: [s]^{k}_a\mapsto a,
$$
where $[s]^{k}_a$ is the $k$-jet of the section $s$ of the bundle $\pi$ at the point $a\in \mathcal{B}$.

The Lie pseudo-group $G$ acts transitively on the base $\mathcal{B}$ of the bundle $\pi$. Because of this without loss of generality we can choose the point $0$ with $y=0, u=0$ on $\mathcal{B}$. Let $N^k$ be the stratum of this point $0$: $N^k=\pi^{-1}(0)$. The stratum $N^k$ is a smooth manifold with coordinates $z, z_y, z_u, z_{yy},z_{yu},\dots, z_{u\dots u}$,
\begin{equation}\label{binom}
\dim N^k={k+2 \choose 2}.
\end{equation}
The stationary subalgebra $\mathcal{G}_0$ of the point $0\in \mathcal{B}$ is generated by the vector fields
$$
Z_1=y\frac{\partial}{\partial y}, \quad  Z_2=z\frac{\partial}{\partial z},\quad   Z_3=u l(u)\frac{\partial}{\partial u},
$$
where $l=l(u)$ is a smooth function. Evolutionary parts (see \cite{VKL1986}) of the vector fields $Z_1^{(k)}, Z_2^{(k)}, Z_3^{(k)}$ are tangent to $N^k$. Let $\overline{Z_1}^{(k)}, \overline{Z_2}^{(k)}, \overline{Z_3}^{(k)}$ be restrictions of these fields to $N^k$. Simple calculations show that in a regular point $\theta\in N^k$ the rank of the system of the tangent vectors $\overline{Z_{1,\theta}}^{(k)}, \overline{Z_{2,\theta}}^{(k)}, \overline{Z_{3,\theta}}^{(k)}$ is equal to $k+2$. Let $G_0$ be the Lie pseudo-group that is generated by subalgebra $\mathcal{G}_0$ and $G_0^{(k)}$ be its prolongation to the space of $k$-jets $J^k(\pi)$. The codimension of the $G_0^{(k)}$-orbit of the point $\theta$ is equal to:
$$
{k+2 \choose 2}-(k+2)=\frac{k^2+k-2}{2}.
$$
But the number of basic differential invariants of order $\leq k$ is a codimension of a regular orbit. So, we obtain formula (\ref{nu}).
\end{pf}

The following theorem gives a complete description of the structure of Petrov differential invariants of QHOE.

 \begin{thm}\label{thm7}

The field of Petrov invariants of QHOE is generated by two second-order differential invariants  $J_{21}$, $J_{22}$ and
invariant  derivations $\nabla_1$ and $\nabla_2$. This field separates regular orbits.
 \end{thm}

\begin{pf}
By lemma \ref{LemmaNu}, the number of basic Petrov invariants of order $k$ ($k\geq 2$) is equal to
$$
\mu(k)=\nu(k)-\nu(k-1)=k.
$$

So, for $k=2$ we have two invariants: $J_{21}$ and $J_{22}$. Applying to them invariant derivations $\nabla_1$ and $\nabla_2$ we obtain four invariants of third order. Since $\mu(3)=3$, there exists one syzygy between them. Indeed, this syzygy we obtain from  (\ref{NablaOnInvariants}):
\begin{equation}\label{syzygy}
 \nabla_2( J_{21})-J_{21} + J_{21} J_{22} = \nabla_1 (J_{22})- J_{22}  + J_{22}^2.  \\
\end{equation}

For $k=4$ we obtain $\mu(4)=4$, but applying to the invariants of third order we can obtain six invariants of fourth order. Therefore there exists two syzygies. We obtain them from (\ref{syzygy}) applying invariant derivations.

So, in order to construct all invariants of order $k+1$ it is sufficient to apply operators $\nabla_1$ and $\nabla_2$ to constructed invariants of $k$-th order.

The fact that Petrov invariants separate regular orbits follows from \cite{Kruglikov}.

\end{pf}

\section{RESULTS AND DISCUSSION}
\par We have  completely described the field of rational differential invariants of one class of second order ordinary differential equations with scalar control parameter with respect to Lie pseudo-group of feedback transformations. In particular, the explicit expressions for differential invariants of different orders and for invariant derivations have been obtained. General expression for number of basic differential invariants up to fifth order has been presented. Normal forms of QHOE have been specified. All the results obtained are local.
\par However, we understand that we have described a limited class of conservative systems. Non-conservative systems are of the great interest. In particular, introducing dissipation to the system can make it rather more interesting for engineering applications, but on same time may lead to a significant complication of its complete description. Furthermore, different types of stochastic terms can be included in the equation considered.  In this case, the powerful apparatus of mathematical statistics can be applied simultaneously with the jet space theory. Finally, the developed approach is applicable for classical systems. However, we believe that possible generalization of the proposed methods may contribute to the field of quantum mechanics and quantum field theory.

\section{APPENDIX}
\subsection{A. Fourth- and fifth-order differential invariants}

\par  Functions
 \begin{equation*}\label{inv4}
\begin{aligned}
J_{41} =&\frac{z_{yyyy}z^3}{z_y^4}, \quad J_{42} =\frac{z_{yyyu}z^3}{z_y^3z_u}, \\
  J_{43} =&\frac{z_{yyuu}z^3}{z_u^2z_y^2}-J_{32}\frac{z_{uu}z}{z_{u}^{2}}, \quad J_{44}=&\frac{z_{yuuu}z^3}{z_y z_u^3}-3J_{33}\frac{z_{uu}z}{z_{u}^{2}} - J_{22}\frac{z_{uuu}z^2}{z_{u}^{3}},
 \end{aligned}
 \end{equation*}
and
\begin{equation*}
\begin{aligned}
J_{51} =&\frac{z_{yyyyy}z^4}{z_y^5}, \quad J_{52}=\frac{z_{yyyyu}z^4}{z_y^4z_u}, \quad J_{53}=\frac{z_{yyyuu}z^4}{z_u^3z_y^2}-J_{42}\frac{z_{uu}z}{z_{u}^{2}}, \\
J_{54}=&\frac{z_{yyuuu}z^4}{z_y^2 z_u^3}-3J_{43}\frac{z_{uu}z}{z_{u}^{2}} - J_{32}\frac{z_{uuu}z^2}{z_{u}^{3}},\\ \label{inv5}
J_{55}=&\frac{z_{yuuuu}z^4}{z_y z_u^4}-3J_{44}\frac{z_{uu}z}{z_{u}^{2}} - J_{33}\Bigl(\frac{z_{uu}z}{z_{u}^{2}}\Bigr)^2-J_{33} \frac{z_{uuu}z^2}{z_{u}^{3}} - J_{22} \frac{z_{uuuu}z^3}{z_{u}^{4}}
 \end{aligned}
\end{equation*}
form a complete set of the basic  fourth- and fifth-order differential invariants,  respectively.

\subsection{B. Application of invariant derivations to third- and fourth-order differential invariants}
\par As we apply invariant derivations $\nabla_{1,2}$ to third- and fourth-order differential invariants following relations (syzygies) are satisfied:
 \begin{align*}
 \nabla_1( J_{31}) & = 2J_{31} - 3 J_{21}J_{31} + J_{41} ,\\
 \nabla_2( J_{31}) & = 2 J_{31} - 3 J_{22} J_{31} + J_{42}, \\
 \nabla_1 (J_{32}) & =  2J_{32} - 2J_{21} J_{32}  - J_{22}J_{32} + J_{42}, \\
 \nabla_2 (J_{32}) & =  2J_{32} - 2J_{22}J_{32} + J_{43}, \\
\nabla_1 (J_{33}) & =  2J_{33} - 2J_{21}J_{33}-3J_{22}J_{33} + J_{43}, \\
 \nabla_2 (J_{33}) & =  2J_{33} - J_{22}J_{33} + J_{44},
 \end{align*}

 \begin{align*}
 \nabla_1( J_{41}) & = 3J_{41} - 4 J_{21}J_{41} + J_{51}, \\
 \nabla_2( J_{41}) & = 3 J_{41} - 4 J_{22} J_{41} + J_{52}, \\
 \nabla_1 (J_{42}) & =  3J_{42} - 3J_{21} J_{42}  - J_{22}J_{42} + J_{52}, \\
\nabla_2 (J_{42}) & =  3J_{42} - 3J_{22}J_{42} + J_{53}, \\
 \nabla_1 (J_{43}) & =  3J_{43} - 2J_{21}J_{43}-2J_{22}J_{43}-J_{33}J_{32} + J_{53}, \\
 \nabla_2 (J_{43}) & =  3J_{43} - 2J_{22}J_{43} + J_{54}, \\
 \nabla_1 (J_{44}) & =  3J_{44} - J_{21}J_{44}-4J_{22}J_{44}-J_{33}^2 + J_{54}, \\
 \nabla_2 (J_{44}) & =  3J_{44} - J_{22}J_{44} + J_{55}.
 \end{align*} 
\end{document}